\documentclass{amsart}

\usepackage[all]{xy}
\usepackage{amscd}
\usepackage[mathscr]{eucal}
\usepackage{amssymb}

\newcommand{\colim}{\operatornamewithlimits{colim}}
\newcommand{\holim}{\operatornamewithlimits{holim}}
\newcommand{\hocolim}{\operatornamewithlimits{hocolim}}
\DeclareMathOperator{\Map}{Map}
\newtheorem{lem}{Lemma}[section]

\newtheorem{prop}[lem]{Proposition}
\newtheorem{theo}[lem]{Theorem}
\newtheorem{coro}[lem]{Corollary}

\newtheorem*{theorem}{Theorem}

\newtheorem*{rema}{Remark}

\begin{document}
\title{Continuous Control and the Algebraic $L$-theory Assembly Map}
\author{David Rosenthal}
\address{Department of Mathematics and Statistics, McMaster University, Hamilton, ON L8S 4K1, Canada}
\email{rosend@math.mcmaster.ca}

\begin{abstract}
	In this work, the assembly map in $L$-theory for the family of finite subgroups is proven to be a split injection for a class of groups. Groups in this class, including virtually polycyclic groups, have universal spaces that satisfy certain geometric conditions. The proof follows the method developed by Carlsson-Pedersen to split the assembly map in the case of torsion free groups. Here, the continuously controlled techniques and results are extended to handle groups with torsion.
\end{abstract}

\maketitle

\section{Introduction}

	The algebraic $L$-groups, $L_n(\mathbb{Z}\Gamma)$, contain geometric information about smooth manifolds whose fundamental group is $\Gamma$. In an effort to calculate these groups, various assembly maps are studied that relate the $L$-groups to a generalized homology theory evaluated on the universal space for a chosen family of subgroups. In this paper, the assembly map for the family of finite subgroups is analyzed. Homology groups are somewhat easier to calculate since tools such as spectral sequences and Chern characters are available. If the assembly map is a {\it split injection}, then the homology groups are direct summands of the $L$-groups. Thus, a piece of the puzzle can be understood. Using {\it continuously controlled algebra}, Carlsson and Pedersen proved that the assembly maps in algebraic $K$- and $L$-theory are split injective for torsion free groups satisfying certain geometric conditions listed below~\cite{cp}. In~\cite{me2}, their $K$-theory result was extended to torsion free groups. The goal here is to extend the $L$-theory result.	
	
	Let $\Gamma$ be a discrete group, and let $\mathcal{E}=E\Gamma(\mathfrak{f})$ be the universal space for $\Gamma$-actions with finite isotropy, where $\mathfrak{f}$ denotes the family of finite subgroups of $\Gamma$. Assume that $\mathcal{E}$ is a finite $\Gamma$-CW complex admitting a compactification, $X$, (i.e., $X$ is compact and $\mathcal{E}$ is an open dense subset) such that
	\begin{itemize}
		\item the $\Gamma$-action extends to $X$;
		\item $X$ is metrizable;
		\item $X^G$ is contractible for every $G\in \mathfrak{f}$;
		\item ${\mathcal{E}}^G$ is dense in $X^G$ for every $G\in \mathfrak{f}$;
		\item compact subsets of $\mathcal{E}$ become small near $Y=X-\mathcal{E}$.  That is, for every compact subset $K \subset \mathcal{E}$ and for every neighborhood $U \subset X$ of $y \in Y$, there exists a neighborhood $V \subset X$ of $y$ such that $g \in \Gamma$ and $gK \cap V \neq \emptyset$ implies $gK \subset U$.
	\end{itemize}
These assumptions are satisfied by a large class of groups, including virtually polycyclic groups and word hyperbolic groups. The following theorem is proven in this paper.

\begin{theorem}
	If $\Gamma$ satisfies the above conditions, and $R$ is a ring with involution such that for sufficiently large $i$, $K_{-i}(RH)=0$ for every $H\in \mathfrak{f}$, then the assembly map, $H_*^{\Gamma} (\mathcal{E};\mathbb{L}^{-\infty}(R{\Gamma}_x)) \to L^{\langle -\infty \rangle}_*(R\Gamma)$, is a split injection.
\end{theorem}

	There are many groups that satisfy these conditions, including virtually polycyclic groups. The reason for the restriction on the ring is explained below. This assumption is not too restrictive since the most important example, in terms of geometric applications, is $R=\mathbb{Z}$, and for $i\geq2$, $K_{-i}(\mathbb{Z}H)=0$ for every finite group $H$.

\subsection{A Brief Outline}

	The idea here is to take the theory of continuously controlled algebra developed to split the algebraic $K$-theory assembly map~\cite{cp, me2} and use it to split the assembly map in algebraic $L$-theory. This was done successfully in the case of torsion free groups in~\cite{cp}. There are two main differences between the $K$- and $L$-theory splittings. The most obvious is that the $\mathbb{K}^{-\infty}$ functor is used in algebraic $K$-theory, and the $\mathbb{L}^{-\infty}$ functor is used in algebraic $L$-theory. The key properties of the $\mathbb{K}^{-\infty}$ functor used in~\cite{cp, me2} are that it takes Karoubi filtrations to fibrations of spectra, and that it commutes with  both products and taking fixed sets. The $\mathbb{L}^{-\infty}$ functor also takes Karoubi filtrations to fibrations of spectra and commutes with taking fixed sets, but it does not behave as well with respect to products. A sufficient condition for $\mathbb{L}^{-\infty}$ to commute with products is that for sufficiently large $i$, $K_{-i}(RH)=0$ for every $H\in \mathfrak{f}$. The other important difference between the two splittings is that for algebraic $L$-theory, rings with involution are used. This induces an involution on the continuously controlled categories. Therefore in order to use the results from the algebraic $K$-theory splitting, we need to make sure that their proofs preserve the involution. As mentioned above, this was verified in~\cite{cp}. In this paper, we wish to do the same with respect to the results in~\cite{me2}.
	
	Section~\ref{l} contains a brief discussion of the $\mathbb{L}^{-\infty}$ functor. The continuously controlled categories are introduced in Section~\ref{control}. In Section~\ref{assembly}, the assembly map is formulated as a map of fixed sets. Realizing the assembly map as a map of fixed sets allows us to use homotopy fixed sets to construct a splitting on the level of homotopy groups. This is explained in Section~\ref{fixed}. The proof of the main theorem will reduce to working with finite group actions, at which point it will be necessary to show that the reduced Steenrod homology (from Theorem~\ref{sthom}) of the quotient of a contractible compact metrizable space by a finite group is trivial. This is discussed in Section~\ref{conner} making use of the Conner Conjecture. Finally, Sections~\ref{proof} and~\ref{filtering} contain the proof of the main theorem using a filtration of the space by {\it conjugacy classes of fixed sets}.

\section{Preliminaries}

\subsection{The $\mathbb{L}^{-\infty}$ Functor}\label{l}

	The functor $\mathbb{L}^{-\infty}$ assigns to an additive category with involution $\mathcal{A}$ (as defined by Ranicki~\cite{rana}) a spectrum whose homotopy groups are the {\it ultimate lower quadratic $L$-groups}, $L^{\langle -\infty \rangle}_*(\mathcal{A})$~\cite[Chapter 17]{ranl}.  There are three properties that this functor possesses (verified in~\cite{cp}) that play an important role in this work.
	
	1. If a group acts on $\mathcal{A}$, then taking fixed sets commutes with applying $\mathbb{L}^{-\infty}$.
	
	2. If $\{ \mathcal{A}_j \}$ is a countable collection of additive categories with involution, then $\mathbb{L}^{-\infty}(\prod \mathcal{A}_j) \simeq \prod \mathbb{L}^{-\infty}(\mathcal{A}_j)$, provided there is an $i_0$ independent of $j$, such that for every $j$, $K_{-i}(\mathcal{A}_j)=0$ whenever $i \geq i_0$.

	3. Let $\mathcal{U}$ be an $\mathcal{A}$-{\it filtered} additive category (as defined by Karoubi~\cite{karoubi}), where $\mathcal{A}$ is a full subcategory of $\mathcal{U}$, and let $\mathcal{U} / \mathcal{A}$ be the associated {\it quotient category}, whose objects are the same as those in $\mathcal{U}$, but whose morphisms are identified if their difference factors through $\mathcal{A}$. If there is an involution on $\mathcal{U}$ that preserves $\mathcal{A}$, then $\mathbb{L}^{-\infty}(\mathcal{A}) \to \mathbb{L}^{-\infty}(\mathcal{U}) \to \mathbb{L}^{-\infty}(\mathcal{U} / \mathcal{A})$ is a homotopy fibration of spectra.

\subsection{Continuously Controlled Algebra}\label{control}

	Let $\Gamma$ be a discrete group, $X$ a $\Gamma$-space, and $Y$ a closed $\Gamma$-invariant subspace of $X$. Let $E=X-Y$ and $R$ be a ring with unit. The free $R$-module generated by $E \times \Gamma \times \mathbb{N}$ is denoted $R[E \times \Gamma]^{\infty}$. If $A$ is a submodule of $R[E \times \Gamma]^{\infty}$, denote $A \cap R[x \times \Gamma]^{\infty}$ by $A_x$, where $x \in E$. The {\it continuously controlled category}, $\mathcal{B}(X,Y;R)$, has objects, $A$, where
\begin{enumerate}
	\item[(i)] $A=\bigoplus_{x \in E} A_x$;
	\item[(ii)] $A_x$ is a finitely generated free $R$-module with basis contained in $\{ x \} \times \Gamma \times \mathbb{N}$;
	\item[(iii)] $\{x \in E \, | \, A_x \neq 0 \}$ is locally finite in $E$.
\end{enumerate}
Morphisms are all $R$-module morphisms, $\phi:A\to B$, that are {\it continuously controlled}. This means that for every $y\in Y$ and every neighborhood $U \subseteq X$ of $y$, there exists a neighborhood $V \subseteq X$ of $y$ such that the components of $\phi$, $\phi^x_z:A_x \to B_z$ and $\phi^z_x: B_z \to A_x$, are zero whenever $x \in V$ and $z \notin U$.

	If $R$ is a ring with involution, then there is an involution, *, on $\mathcal{B}(X,Y;R)$ defined as follows. For every object $A$, $A^*=A$. The component $(\phi^*)^x_z:B^*_x \to A^*_z$ of the morphism $\phi^*:B^* \to A^*$, is defined to be the conjugate transpose of $\phi^z_x:A_z \to B_x$, where the elements of $R$ are conjugated by the involution on the ring. Here $A_x$ is identified with its dual via its basis, and the matrix description of the dual map is used.

	The diagonal action of $\Gamma$ on $E \times \Gamma$ induces a $\Gamma$-action on $R[E \times \Gamma]^{\infty}$ that induces a $\Gamma$-action on $\mathcal{B}(X,Y;R)$. If we assume that $\Gamma$ acts on $E$ with finite isotropy, then the fixed category, $\mathcal{B}^{\Gamma}(X,Y;R)$, has those objects, $A$, in $\mathcal{B}(X,Y;R)$ that satisfy the conditions
\begin{enumerate}
  \item[1.] $A_{\gamma x}\cong A_x$ for every $\gamma \in \Gamma$, and
  \item[2.] $A_x$ is a finitely generated free $R{\Gamma}_x$-module. 
\end{enumerate}  
This implies that $\bigoplus_{x'\in [x]} A_{x'}$ is a finitely generated free $R\Gamma$-module, where $[x]=\{ \gamma x \, | \, \gamma \in \Gamma \}$. The morphisms in $\mathcal{B}^{\Gamma}(X,Y;R)$ are those morphisms, $\phi$, in $\mathcal{B}(X,Y;R)$ that are $\Gamma$-equivariant, i.e., $\gamma {\phi}^x_y {\gamma}^{-1}=\phi^{\gamma x}_{\gamma y}$ for all $\gamma \in \Gamma$ and all $x,y \in E$. Notice that if there are infinite isotropy subgroups, then every object in the fixed category will be the zero object.

\begin{rema}
{\rm		In~\cite{hp}, the {\it equivariant} continuously controlled category, $\mathcal{B}_{\Gamma}(X,Y;R)$, is defined. This category has a slightly different definition from the one used here. These categories yield equivariant homology theories which are used to formulate assembly maps for any given family of subgroups. In this work, we are only concerned with the case of the family of finite subgroups and so we will always assume that $\Gamma$ acts on $E$ with finite isotropy. With this assumption, the equivariant category, $\mathcal{B}_{\Gamma}(X,Y;R)$, is equivalent to the fixed category, $\mathcal{B}^{\Gamma}(X,Y;R)$, which makes it possible to realize the assembly map as a map of fixed sets. Doing so provides a way to construct a splitting using homotopy fixed sets (see Section~\ref{fixed}). In this paper, we will use the subscript ${-}_{\Gamma}$, as in $\mathcal{B}_{\Gamma}(X,Y;R)$, for the equivariant continuously controlled categories, and the superscript ${-}^{\Gamma}$, as in $\mathcal{B}^{\Gamma}(X,Y;R)$, for the fixed categories described above. }
\end{rema}
	
	Two controlled categories that are of particular interest to us are $\mathcal{B}(CX,CY \cup X,p_X;R)$ and $\mathcal{B}(\Sigma X,\Sigma Y,p_X;R)$, where $CX$ denotes the cone of $X$, $\Sigma X$ denotes the unreduced suspension of $X$, and $p_X: X\times (0,1) \to X$ is the projection map. These categories have the same objects as $\mathcal{B}(CX,CY \cup X;R)$ and $\mathcal{B}(\Sigma X,\Sigma Y;R)$ respectively, but their control conditions on morphisms differ along $Y \times (0,1)$, where they are only required to be $p_X$-{\it controlled}. This means that for every $(y,t)\in Y \times (0,1)$ and every neighborhood $U \subseteq X$ of $y$, there is a neighborhood $V$ of $(y,t)$ such that $\phi^a_b =0$ and $\phi^b_a =0$ whenever $a \in V \cap p_X^{-1}(U)$ and $b \notin p_X^{-1}(U)$.

	The {\it support at infinity} of an object $A$ in $\mathcal{B}(X,Y;R)$, denoted ${\rm supp}_{\infty}(A)$, is the set of limit points of $\{x \, | \, A_x \neq 0\}$. If $C$ is a closed  subspace of $Y$, then the category $\mathcal{B}(X,Y;R)_C$ is the full subcategory of $\mathcal{B}(X,Y;R)$ on objects, $A$, with ${\rm supp}_{\infty}(A)\subseteq C$. If $W$ is an open subspace of $Y$, then the {\it germ category}, $\mathcal{B}(X,Y;R)^W$, has the same objects as $\mathcal{B}(X,Y;R)$, but morphisms are identified if they agree in a neighborhood of $W$. It is equivalent to the quotient category $\mathcal{B}(X,Y;R) / \mathcal{B}(X,Y;R)_C$, when $C=Y-W$~\cite{cp}. Thus
\[ \mathcal{B}(X,Y;R)_C \to \mathcal{B}(X,Y;R) \to \mathcal{B}(X,Y;R)^{Y-C} \]
induces a fibration of spectra after applying $\mathbb{L}^{-\infty}$. Such sequences therefore become the main tool in this theory.

	A function $f:(X,Y) \to (X',Y')$ is {\it eventually continuous} if for every compact $K \subset X'-Y'$, $f^{-1}(K)$ has compact closure in $X-Y$, $f(X-Y) \subset X'-Y'$, and $f$ is continuous on $Y$. Such a map induces a functor $f_{\sharp}:\mathcal{B}(X,Y;R) \to \mathcal{B}(X',Y';R)$, up to natural equivalence, defined as follows. An object $A$ is sent to $f_{\sharp}A$, where $(f_{\sharp}A)_{x'}=\bigoplus^{}_{x\in f^{-1}(x')}A_x$, and morphisms are induced by the identity. If $f$ sends $Y-W$ to $Y'-W'$, then $f$ also induces a functor $f_{\sharp}:\mathcal{B}(X,Y;R)^W \to \mathcal{B}(X',Y';R)^{W'}$. Note that if $R$ is a ring with involution and $\phi$ is a morphism in $\mathcal{B}(X,Y;R)$, then $f_{\sharp}(\phi^*)=f_{\sharp}(\phi)^*$. Furthermore, the proof of~\cite[Proposition 1.18]{cp} is easily seen to preserve the involution which says that if $f_1$, $f_2:(X,Y) \to (X',Y')$ are eventually continuous maps and $f_1|_Y=f_2|_Y$, then $f_{1\sharp}$ and $f_{2\sharp}$ are naturally equivalent functors. This simple observation is very useful when extending those results from the $K$-theory case that rely on eventually continuous maps. The above discussion also works for fixed categories if we use $\Gamma$-invariant subspaces and $\Gamma$-equivariant functions.

\subsection{The Continuously Controlled Assembly Map}\label{assembly}

	Let $R$ be a ring with involution. Since the Karoubi filtrations described in the previous section induce fibrations of spectra after applying $\mathbb{L}^{-\infty}$, one can prove that the functor
\[ \mathbb{L}^{-\infty}(\mathcal{B}_{\Gamma}(-\times [0,1],- \times 1;R)^{-\times 1}) \]
from the category of $\Gamma$-CW complexes to the category of spectra, is $\Gamma$-excisive and $\Gamma$-homotopy invariant~\cite{hp}. This implies that it satisfies the Eilenberg-Steenrod axioms, except for the dimension axiom, and therefore defines a generalized homology theory. We define
\[ H_*^{\Gamma} (E;\mathbb{L}^{-\infty}(R{\Gamma}_x)) = \pi_*(\Omega\mathbb{L}^{-\infty}(\mathcal{B}_{\Gamma}(E\times [0,1],E\times 1;R)^{E\times 1})). \]
This homology theory is isomorphic to the homology of $E$ over the orbit category~\cite{hp}, defined by Davis and L\"{u}ck.

	The continuously controlled assembly map is defined in~\cite{hp} to be the map of spectra
\[ \Omega\mathbb{L}^{-\infty}(\mathcal{B}_{\Gamma}(E\times [0,1],E\times 1;R)^{E\times 1}) \to \Omega\mathbb{L}^{-\infty}(\mathcal{B}_{\Gamma}(\bullet \times [0,1],\bullet \times 1;R)^{\bullet \times 1}), \]
induced by collapsing $E$ to a point. Hambleton and Pedersen proved that the continuously controlled assembly map is homotopy equivalent to the Davis-L\"{u}ck assembly map~\cite[Theorem 8.3]{hp}.

	Let $\mathcal{F}$ be a family of subgroups of $\Gamma$ that is closed under conjugation and under the operation of taking subgroups. Then there is a universal space $E\Gamma(\mathcal{F})$ for $\Gamma$-actions with isotropy in $\mathcal{F}$. The space $E\Gamma(\mathcal{F})$ is a $\Gamma$-CW complex that is characterized by the fact that $E\Gamma(\mathcal{F})^G$ is contractible for all $G\in \mathcal{F}$ and is empty otherwise. In this paper, the assembly map of interest is the one induced by collapsing $E\Gamma(\mathfrak{f})$ to a point, where $\mathfrak{f}$ denotes the family of finite subgroups. Our goal is to prove that it is a split injection on homotopy groups. The method for doing this follows~\cite{cp}. In this approach, it is necessary to realize the assembly map as a map of fixed sets. To do so, a few assumptions need to be placed on $E\Gamma(\mathfrak{f})$.
	
	Assume that $\Gamma$ acts cocompactly on $\mathcal{E}=E\Gamma(\mathfrak{f})$, and that $\mathcal{E}$ admits a metrizable compactification $X$, such that the $\Gamma$-action extends to $X$, and compact subsets of $\mathcal{E}$ become small near $Y=X-\mathcal{E}$ (i.e., for every compact subset $K \subset \mathcal{E}$ and for every neighborhood $U \subset X$ of $y \in Y$, there exists a neighborhood $V \subset X$ of $y$ such that $g \in \Gamma$ and $gK \cap V \neq \emptyset$ implies $gK \subset U$). Then using an argument as in~\cite[Lemma 2.3]{cp}, one proves that $\Omega\mathbb{L}^{-\infty}(\mathcal{B}^{\Gamma}(\Sigma X,\Sigma Y,p_X;R))$ is weakly homotopy equivalent to $\mathbb{L}^{-\infty}(R\Gamma)$. Furthermore, the map
\[ \mathbb{L}^{-\infty}(\mathcal{B}^{\Gamma}(CX,CY \cup X,p_X;R)) \to \mathbb{L}^{-\infty}(\mathcal{B}^{\Gamma}(\Sigma X,\Sigma Y,p_X;R)), \]
induced by collapsing $X$ from $CX$ to $\Sigma X$, is weakly homotopy equivalent to the map $\mathbb{L}^{-\infty}(\mathcal{B}_{\Gamma}(\mathcal{E}\times [0,1],\mathcal{E}\times 1;R)^{\mathcal{E}\times 1}) \to \mathbb{L}^{-\infty}(\mathcal{B}_{\Gamma}(\bullet \times [0,1],\bullet \times 1;R)^{\bullet \times 1})$~\cite{cp}. Since $ \mathbb{L}^{-\infty}$ commutes with fixed sets, the assembly map can now be realized as the map of fixed spectra
\[ \Omega\mathbb{L}^{-\infty}(\mathcal{B}(CX,CY \cup X,p_X;R))^{\Gamma} \to \Omega\mathbb{L}^{-\infty}(\mathcal{B}(\Sigma X,\Sigma Y,p_X;R))^{\Gamma}. \]

\subsection{Generalized Fixed Sets}\label{fixed}

	In this section, we recall some facts about generalized homotopy fixed sets in the category of spectra, and illustrate how they help us exploit the fact that the assembly map is a map of fixed sets. Let $S$ be a spectrum with $\Gamma$-action. The fixed set $S^{\Gamma}$ can be identified with the set $\Map_{\Gamma}(\bullet,S)$ of $\Gamma$-equivariant maps from a point into $S$. Let $\mathcal{F}$ be a family of subgroups of $\Gamma$. The homotopy fixed set associated to this family, $S^{h \mathcal{F} \Gamma}$, is defined to be $\Map_{\Gamma}(E\Gamma(\mathcal{F}),S)$. Notice that the $\Gamma$-equivariant map $E\Gamma(\mathcal{F}) \to \bullet$
induces the following commutative diagram:
\[ \xymatrix{
	S^{\Gamma} \ar[d]_a \ar[r] & T^{\Gamma} \ar[d] \\
	S^{h\mathcal{F} \Gamma} \ar[r]^b & T^{h\mathcal{F} \Gamma}. } \]
If we can show that $a$ and $b$ are weak homotopy equivalences, then on the level of homotopy groups, we will have successfully split the map $S^{\Gamma} \to T^{\Gamma}$. The next two lemmas play an important role in achieving this goal. The author would like to thank Wolfgang L\"{u}ck for showing him the elegant proof of the following lemma.

\begin{lem}\label{lemmaa}
	Let $F:S \to T$ be an equivariant map of spectra with $\Gamma$-action. If $F^G:S^G \to T^G$ is a weak homotopy equivalence for every $G\in \mathcal{F}$, then $S^{h \mathcal{F} \Gamma} \simeq T^{h \mathcal{F} \Gamma}$.
\end{lem}

\begin{proof}
	Suppose for the moment that $\Map_{\Gamma}(X,S) \to \Map_{\Gamma}(X,T)$ is a weak homotopy equivalence for any finite dimensional $\Gamma$-CW complex $X$. Since the inclusion of the $n$-skeleton of $E\Gamma(\mathcal{F})$ into its $(n+1)$-skeleton $E\Gamma(\mathcal{F})_n \hookrightarrow E\Gamma(\mathcal{F})_{n+1}$, is a cofibration, $E\Gamma(\mathcal{F})=\colim_n E\Gamma(\mathcal{F})_n \simeq \hocolim_n E\Gamma(\mathcal{F})_n$. Therefore
\begin{align}
	&S^{h \mathcal{F} \Gamma} & &\simeq & & \Map_{\Gamma}\big(\hocolim_n E\Gamma(\mathcal{F})_n,S\big) & &\simeq & &\holim_n\Map_{\Gamma}(E\Gamma(\mathcal{F})_n,S) \notag\\
	& & & & & & &\simeq & &\holim_n\Map_{\Gamma}(E\Gamma(\mathcal{F})_n,T)\notag\\
	& & & & & & &\simeq & &T^{h \mathcal{F} \Gamma}.\notag
\end{align}	

	Let $X$ be a finite dimensional $\Gamma$-CW complex. We will prove that $\Map_{\Gamma}(X,S) \simeq \Map_{\Gamma}(X,T)$ by induction on the skeleta of $X$. 
\[ \Map_{\Gamma}(X_0,S)=\Map_{\Gamma}\Big(\bigsqcup_i \Gamma / H_i,S\Big)=\prod_i\Map_{\Gamma}(\Gamma / H_i,S)=\prod_i S^{H_i}, \]
where each $H_i \in \mathcal{F}$. Since $S^{H_i}\simeq T^{H_i}$ for every $i$, $\Map_{\Gamma}(X_0,S) \simeq \Map_{\Gamma}(X_0,T)$.

	Now assume that $\Map_{\Gamma}(X_{k-1},S) \simeq \Map_{\Gamma}(X_{k-1},T)$. Consider the following pushout diagram
\[ \xymatrix{
	\Map_{\Gamma}\big(\bigsqcup_j \Gamma / G_j \times S^{k-1},S\big) \ar[d] \ar[r] & \Map_{\Gamma}(X_{k-1},S) \ar[d] \\
	\Map_{\Gamma}\big(\bigsqcup_j \Gamma / G_j \times \mathbb{D}^k,S\big) \ar[r] & \Map_{\Gamma}(X_k,S) } \]
where each $G_j\in \mathcal{F}$. The equivariant map from $S$ to $T$ induces a commutative diagram relating the above pushout diagram to the corresponding one for $T$. Therefore the proof will be complete if we show that the maps $\Map_{\Gamma}(X_{k-1},S) \to \Map_{\Gamma}(X_{k-1},T)$ and $\Map_{\Gamma}\big(\bigsqcup_j \Gamma / G_j \times \mathbb{D}^k,S\big) \to \Map_{\Gamma}\big(\bigsqcup_j \Gamma / G_j \times \mathbb{D}^k,T\big)$ are weak homotopy equivalences. The first map is a weak homotopy equivalence by the induction hypothesis, and since
\[ \Map_{\Gamma}\big(\bigsqcup_j \Gamma / G_j \times \mathbb{D}^k,S\big)=\prod_j\Map_{\Gamma}(\Gamma / G_j \times \mathbb{D}^k,S) \simeq \prod_j\Map_{\Gamma}(\Gamma / G_j,S)= \prod_j S^{G_j} \]
the second map is also a weak homotopy equivalence.
\end{proof}

\begin{lem}\label{lemmab}
	Let $B$ be a $G$-spectrum, where $G\in \mathcal{F}$, and let $\Gamma$ act on $S=\prod_{\Gamma / G}B$ by identifying $S$ with $\Map_G(\Gamma,B)$ (in which $(\gamma f)(x)=f({\gamma}^{-1}x)$, where $\gamma \in \Gamma$). Then $S^{\Gamma} \simeq S^{h \mathcal{F} \Gamma}$.
\end{lem}

\begin{proof}
\[ S^{h \mathcal{F} \Gamma}= \Map_{\Gamma}(E\Gamma(\mathcal{F}),\Map_G(\Gamma,B))=\Map_G(E\Gamma(\mathcal{F}),B) \]
by evaluating at 1.
Since $E\Gamma(\mathcal{F})$ is $G$-equivariantly homotopy equivalent to $E\Gamma(\mathcal{F})^G$ and $E\Gamma(\mathcal{F})^G$ is contractible,
\begin{align}
	&\Map_G(E\Gamma(\mathcal{F}),B) & &\cong & &\Map_G(E\Gamma(\mathcal{F})^G,B) & &= & &\Map(E\Gamma(\mathcal{F})^G,B^G) \notag\\
	& & & & & & &\cong & &B^G\notag\\
	& & & & & & &= & &S^{\Gamma}.\notag
\end{align}
\end{proof}

\subsection{The Conner Conjecture and Steenrod Homology}\label{conner}

	Recall that a {\it reduced Steenrod homology theory}, $h$, is a functor from the category of compact metrizable spaces and continuous maps to the category of graded abelian groups satisfying
\begin{enumerate}
	\item[(i)] $h$ is homotopy invariant;
	\item[(ii)] $h(\bullet)=0$;
	\item[(iii)] given any closed subset $A$ of $X$, there is a natural transformation, $\partial :h_n(X / A) \to h_{n-1}(A)$, fitting into a long exact sequence
\[ \dotsb \to h_n(A) \to h_n(X) \to h_n(X / A) \to h_{n-1}(A) \to \dotsb ; \]
	\item[(iv)] if $\bigvee^{}_{} X_i$ denotes a compact metric space that is a countable union of metric spaces along a single point, then the projection maps $p_i : \bigvee^{}_{} X_i \to X_i$ induce an isomorphism $h_*(\bigvee^{}_{} X_i) \to \prod^{}_{} h_*(X_i)$.
\end{enumerate}
Given any generalized homology theory, there is a unique Steenrod homology extension. Existence of such extensions was proved by Kahn, Kaminker and Schochet, and Edwards and Hastings~\cite{kahn,edwards}. Uniqueness was proved by Milnor \cite{milnor}.

	A functor $k$ from the category of compact metrizable spaces to the category of spectra is called a {\it reduced Steenrod functor} if it satisfies the following conditions.
\begin{enumerate}
	\item[(i)] The spectrum $k(CX)$ is contractible.
	\item[(ii)] If $A \subset X$ is closed, then $k(A) \to k(X) \to k(X / A)$ is a fibration (up to natural weak homotopy equivalence).
	\item[(iii)] The projection maps $p_i : \bigvee^{}_{} X_i \to X_i$ induce a weak homotopy equivalence $k(\bigvee^{}_{} X_i) \to \prod^{}_{} k(X_i)$.
\end{enumerate}

\begin{theo}\cite[Theorem 5.4]{cp}\label{sthom}
        Let $R$ be a ring with involution such that $K_{-i}(R)=0$ for sufficiently large $i$. Then the functor \[ \mathbb{L}^{-\infty}(\mathcal{B}(C(-),-;R)) \] is the reduced Steenrod functor associated to the 4-periodic $\mathbb{L}^{-\infty}(R)$-spectrum.
\end{theo}

	The following proposition plays a key role in the proof of Theorem~\ref{filtration}. It follows from \cite[Proof of Theorem 2.13]{cp} and Theorem~\ref{sthom}.

\begin{prop}\label{key}
	Let $X$ be a compact metrizable space, $Y$ a closed nowhere dense subset, and $R$ a ring with involution such that $K_{-i}(R)=0$ for sufficiently large $i$. If the reduced Steenrod homology of $X$ (from Theorem~\ref{sthom}) is trivial, then
\[ \mathbb{L}^{-\infty}(\mathcal{B}(CX,CY \cup X,p_X;R)) \simeq \mathbb{L}^{-\infty}(\mathcal{B}(\Sigma X,\Sigma Y,p_X;R)). \]
\end{prop}

	In the proof of Theorem~\ref{filtration} we will be considering spaces with a finite group action. In order to use Proposition~\ref{key} on the corresponding quotient spaces, we need to know that the reduced Steenrod homology (from Theorem~\ref{sthom}) of the quotient of a contractible compact metrizable space by a finite group is trivial. The first step toward proving this is~\cite[Theorem III.7.12]{bredon}, which states that if $G$ is a finite group acting on a compact Hausdorff space $X$ that has trivial reduced \v{C}ech cohomology, then the reduced \v{C}ech cohomology of $X / G$ is also trivial. This is a special case of the Conner Conjecture which states that if a compact Lie group $G$ acts on a space $X$, where $X$ is either paracompact of finite cohomological dimension with finitely many orbit types or compact Hausdorff, then $X / G$ has trivial reduced \v{C}ech cohomology if $X$ does \cite{conner}. The Conner Conjecture was proved by Oliver~\cite{oliver}. The final step is the following theorem. 

\begin{theo}\cite[Theorem 5.2]{me2}\label{geogh}
	If the reduced \v{C}ech cohomology of a compact metrizable space $X$ is trivial, then  every reduced Steenrod homology of $X$ is also trivial.
\end{theo}

	As mentioned above, the following corollary will be used in the proof of Theorem~\ref{filtration}.

\begin{coro}\label{zero}
	Let $G$ be a finite group acting on a compact metrizable space $X$. If the reduced \v{C}ech cohomology of $X$ is trivial, then the reduced Steenrod homology of $X / G$ (from Theorem~\ref{sthom}) is trivial.
\end{coro}

\section{The Main Theorem}\label{main}

	Let $\Gamma$ be a discrete group, and let $\mathcal{E}=E\Gamma(\mathfrak{f})$ be the universal space for $\Gamma$-actions with finite isotropy, where $\mathfrak{f}$ denotes the family of finite subgroups of $\Gamma$. Assume that $\mathcal{E}$ is a finite $\Gamma$-CW complex admitting a compactification, $X$, (i.e., $X$ is compact and $\mathcal{E}$ is an open dense subset) such that
	\begin{itemize}
		\item the $\Gamma$-action extends to $X$;
		\item $X$ is metrizable;
		\item $X^G$ is contractible for every $G\in \mathfrak{f}$;
		\item ${\mathcal{E}}^G$ is dense in $X^G$ for every $G\in \mathfrak{f}$;
		\item compact subsets of $\mathcal{E}$ become small near $Y=X-\mathcal{E}$.  That is, for every compact subset $K \subset \mathcal{E}$ and for every neighborhood $U \subset X$ of $y \in Y$, there exists a neighborhood $V \subset X$ of $y$ such that $g \in \Gamma$ and $gK \cap V \neq \emptyset$ implies $gK \subset U$.
	\end{itemize}

\begin{theorem}
	If $\Gamma$ satisfies the above conditions, and $R$ is a ring with involution such that for sufficiently large $i$, $K_{-i}(RH)=0$ for every $H\in \mathfrak{f}$, then the assembly map, $H_*^{\Gamma} (\mathcal{E};\mathbb{L}^{-\infty}(R{\Gamma}_x)) \to L^{\langle -\infty \rangle}_*(R\Gamma)$, is a split injection.
\end{theorem}

	There are many important classes of groups satisfying the conditions of the main theorem. {\it Crystallographic groups}, which are discrete groups that act cocompactly on Euclidean $n$-space by isometries, satisfy these conditions. The desired compactification is obtained by adding an $(n-1)$-sphere at infinity. More generally, {\it virtually polycyclic groups} satisfy these conditions, since we can also take $\mathbb{R}^n$ for some $n$ to be our universal space~\cite{wilking}. Gromov's {\it word hyperbolic groups} also satisfy the conditions by taking a certain compactification of the Rips complex~\cite{bestvina, meintrup}.

	For the remainder of this paper, assume that $\Gamma$, $R$, $\mathcal{E}$, $X$, and $Y$ satisfy the conditions of the main theorem.

\subsection{The Proof}\label{proof}

	The assembly map, as shown in Section~\ref{assembly}, is the map of fixed sets
\[ \Omega\mathbb{L}^{-\infty}(\mathcal{B}(CX,CY \cup X,p_X;R))^{\Gamma} \to \Omega\mathbb{L}^{-\infty}(\mathcal{B}(\Sigma X,\Sigma Y,p_X;R))^{\Gamma}. \]
This allows us to use homotopy fixed sets to show that the induced map on homotopy groups is a split injection. We do this by considering the commutative diagram discussed in Section~\ref{fixed}. Then split injectivity of the assembly map is proven by showing that the maps
\[ \Omega\mathbb{L}^{-\infty}(\mathcal{B}(CX,CY \cup X,p_X;R))^{\Gamma} \to \Omega\mathbb{L}^{-\infty}(\mathcal{B}(CX,CY \cup X,p_X;R))^{h\mathfrak{f}\Gamma} \]
and
\[ \Omega\mathbb{L}^{-\infty}(\mathcal{B}(CX,CY \cup X,p_X;R))^{h\mathfrak{f}\Gamma} \to \Omega\mathbb{L}^{-\infty}(\mathcal{B}(\Sigma X,\Sigma Y,p_X;R))^{h\mathfrak{f}\Gamma} \]
are each weak homotopy equivalences. In this section we tackle the first of these two maps as in~\cite[Section 6]{me2}. The other map is handled in Section~\ref{filtering}.

\begin{theo}\label{theo3}
	The spectrum $\Omega\mathbb{L}^{-\infty}(\mathcal{B}(CX,CY \cup X,p_X;R))^{\Gamma}$ is weakly homotopy equivalent to $\Omega\mathbb{L}^{-\infty}(\mathcal{B}(CX,CY \cup X,p_X;R))^{h\mathfrak{f}\Gamma}$.
\end{theo}

\begin{proof}
	It suffices to prove that $\mathbb{L}^{-\infty}(\mathcal{B}(\mathcal{E} \times (0,1],\mathcal{E} \times 1;R)^{\mathcal{E} \times 1})^{\Gamma}$ is weakly homotopy equivalent to $\mathbb{L}^{-\infty}(\mathcal{B}(\mathcal{E} \times (0,1],\mathcal{E} \times 1;R)^{\mathcal{E} \times 1})^{h \mathfrak{f} \Gamma}$. Given this, the theorem is proved as follows:
\begin{align}
	&\mathbb{L}^{-\infty}(\mathcal{B}(CX,CY \cup X,p_X;R))^{\Gamma} & &\simeq & &\mathbb{L}^{-\infty}(\mathcal{B}(CX,CY \cup X,p_X;R)^\mathcal{E})^{\Gamma} \notag\\ 
	& & &= & &\mathbb{L}^{-\infty}(\mathcal{B}(\mathcal{E} \times (0,1],\mathcal{E} \times 1;R)^{\mathcal{E} \times 1})^{\Gamma} \notag\\
	& & &\simeq & &\mathbb{L}^{-\infty}(\mathcal{B}(\mathcal{E} \times (0,1],\mathcal{E} \times 1;R)^{\mathcal{E} \times 1})^{h \mathfrak{f} \Gamma} \notag\\
	& & &= & &\mathbb{L}^{-\infty}(\mathcal{B}(CX,CY \cup X,p_X;R)^\mathcal{E})^{h \mathfrak{f} \Gamma} \notag\\
	& & &\simeq & &\mathbb{L}^{-\infty}(\mathcal{B}(CX,CY \cup X,p_X;R))^{h \mathfrak{f} \Gamma}. \notag
\end{align}
These equivalences follow from~\cite[Lemmas 2.4 and 2.5]{cp}, whose proofs preserve the involution and therefore hold for $L$-theory, and Lemma~\ref{lemmaa}. 

	We now prove that $\mathbb{L}^{-\infty}(\mathcal{B}(\mathcal{E} \times (0,1],\mathcal{E} \times 1;R)^{\mathcal{E} \times 1})^{\Gamma}$ is weakly homotopy equivalent to $\mathbb{L}^{-\infty}(\mathcal{B}(\mathcal{E} \times (0,1],\mathcal{E} \times 1;R)^{\mathcal{E} \times 1})^{h \mathfrak{f} \Gamma}$. Note that $\mathcal{E}$ is assumed to be a finite $\Gamma$-CW complex, and proceed by induction on the $\Gamma$-cells in $\mathcal{E}$. Begin with the discrete space ${\Gamma}/H$ for some $H \in \mathfrak{f}$. The control condition implies that the components of a morphism near ${\Gamma}/H \times 1$ must be zero between points with different ${\Gamma}/H$ entries. Since we are taking germs at ${\Gamma}/H \times 1$, the category $\mathcal{B}({\Gamma}/H \times (0,1],{\Gamma}/H \times 1;R)^{{\Gamma}/H \times 1}$ is equivalent to the product category $\prod_{\Gamma /H}^{} \mathcal{B}((0,1],1;R)^{1}$. The action of $\Gamma$ that the product category inherits is the same as the one on the product in Lemma~\ref{lemmab}. The projection maps induce a map
\[ \mathbb{L}^{-\infty}\Big(\prod_{\Gamma /H}^{} \mathcal{B}((0,1],1;R)^{1}\Big) \to \prod_{\Gamma /H}^{} \mathbb{L}^{-\infty}(\mathcal{B}((0,1],1;R)^{1}) \]
that is $\Gamma$-equivariant, and a weak homotopy equivalence by the assumption that for sufficiently large $i$, $K_{-i}(RH')=0$ for every $H'\in \mathfrak{f}$~\cite{cp}.

	Consider the following commutative diagram:
\[ \xymatrix{
	\mathbb{L}^{-\infty} \Big(\prod_{\Gamma /H}^{} \mathcal{B}((0,1],1;R)^{1}\Big)^{\Gamma} \ar[d]_a \ar[r]^b &  \Big(\prod_{\Gamma /H}^{} \mathbb{L}^{-\infty}(\mathcal{B}((0,1],1;R)^{1})\Big)^{\Gamma} \ar[d]_c \\
	\mathbb{L}^{-\infty} \Big(\prod_{\Gamma /H}^{} \mathcal{B}((0,1],1;R)^{1}\Big)^{h\mathfrak{f} \Gamma} \ar[r]^d & \Big(\prod_{\Gamma /H}^{} \mathbb{L}^{-\infty}(\mathcal{B}((0,1],1;R)^{1})\Big)^{h\mathfrak{f} \Gamma} } \]
We want to show that $a$ is a weak homotopy equivalence.

	Let $G \leq \Gamma$ be given, and choose representatives $\gamma_j$, so that $G \backslash \Gamma /H=\{ G\gamma_jH \}$. Then 
\begin{align}
	&\mathbb{L}^{-\infty} \Big(\prod_{\Gamma /H} \mathcal{B}((0,1],1;R)^{1}\Big)^G & &\cong & &\mathbb{L}^{-\infty} \Big( \Big(\prod_{\Gamma /H}^{} \mathcal{B}((0,1],1;R)^{1} \Big)^G \Big) \notag\\ 
	& & &\cong & &\mathbb{L}^{-\infty} \Big(\prod_{G \backslash \Gamma /H} \mathcal{B}((0,1],1;R[\gamma_j^{-1}G\gamma_j \cap H])^{1} \Big) \notag\\
	& & &\simeq & &\prod_{G \backslash \Gamma /H} \mathbb{L}^{-\infty} (\mathcal{B}((0,1],1;R)^{1})^{\gamma_j^{-1}G\gamma_j \cap H} \notag\\
	& & &\cong & &\Big(\prod_{\Gamma /H} \mathbb{L}^{-\infty} (\mathcal{B}((0,1],1;R)^{1}) \Big)^G \notag
\end{align}
by our assumption that for sufficiently large $i$, $K_{-i}(RH')=0$ for each $H'\in \mathfrak{f}$, and the fact that $\mathbb{L}^{-\infty}$ commutes with taking fixed sets. The case $G=\Gamma$ proves that $b$ is a weak homotopy equivalence. By Lemma $\ref{lemmaa}$, $d$ is a weak homotopy equivalence. Finally, $c$ is a weak homotopy equivalence by Lemma $\ref{lemmab}$.

	Now assume that the theorem holds for $N$ and that $E$ is obtained from $N$ by attaching a $\Gamma$-$n$-cell, $\Gamma / K \times \mathbb{D}^n$, for some $K \in \mathfrak{f}$. Since $\mathcal{B}^{\Gamma}(E \times (0,1],E \times 1;R)_{N \times 1}^{E \times 1}$ is equivalent to $\mathcal{B}^{\Gamma}(N \times (0,1],N \times 1;R)^{N \times 1}$, $\mathcal{B}^{\Gamma}(N \times (0,1],N \times 1;R)^{N \times 1} \to \mathcal{B}^{\Gamma}(E \times (0,1],E \times 1;R)^{E \times 1} \to \mathcal{B}^{\Gamma}(E \times (0,1],E \times 1;R)^{(E-N) \times 1}$ is a Karoubi filtration. Since these are fixed categories and $\mathbb{L}^{-\infty}$ commutes with taking fixed sets, let $A^{\Gamma} \to B^{\Gamma} \to C^{\Gamma}$ denote $\mathbb{L}^{-\infty}$ applied to this sequence. Then $A^{\Gamma} \to B^{\Gamma} \to C^{\Gamma}$ is a fibration of spectra. Since $N$ is $\Gamma$-invariant, it is also $G$-invariant for every $G \in \mathfrak{f}$. Hence, $A^G \to B^G \to C^G$ is a fibration for every $G \in \mathfrak{f}$. Let $D$ be the homotopy fiber of $B \to C$. Taking fixed sets and taking homotopy fibers commute since both are inverse limits. Therefore $D^G \to B^G \to C^G$ is a fibration for every $G \in \mathfrak{f}$. Furthermore $A^G \simeq D^G$ for every $G \in \mathfrak{f}$. Thus, $A^{h \mathfrak{f} \Gamma} \simeq D^{h \mathfrak{f} \Gamma}$, by Lemma $\ref{lemmaa}$. Taking homotopy fixed sets is also an inverse limit. Therefore, it too commutes with taking homotopy fibers. This implies that $D^{h \mathfrak{f} \Gamma} \to B^{h \mathfrak{f} \Gamma} \to C^{h \mathfrak{f} \Gamma}$ is a fibration. As a result, so is $A^{h \mathfrak{f} \Gamma} \to B^{h \mathfrak{f} \Gamma} \to C^{h \mathfrak{f} \Gamma}$.
	
	Consider the following commutative diagram:
\[ \xymatrix{
	A^{\Gamma} \ar[d]_a \ar[r] & B^{\Gamma} \ar[d]_b \ar[r] & C^{\Gamma} \ar[d]_c \\
	A^{h \mathfrak{f} \Gamma} \ar[r] & B^{h \mathfrak{f} \Gamma} \ar[r] & C^{h \mathfrak{f} \Gamma} } \]
We need to show that $b$ is a weak homotopy equivalence. By our induction hypothesis, $a$ is a weak homotopy equivalence. Therefore, by the Five Lemma, it suffices to prove that $c$ is a weak homotopy equivalence. Since $E-N= \Gamma / K \times \mathring{e^n}$, where $\mathring{e^n}$ is an open $n$-cell, the category $\mathcal{B} (E \times (0,1],E \times 1;R)^{(E-N) \times 1}$ is equivalent to the product category $\prod_{\Gamma /K}^{} \mathcal{B} (\mathring{e^n} \times (0,1],\mathring{e^n} \times 1;R)^{\mathring{e^n} \times 1}$. But this is entirely similar to the start of the induction. This completes the proof.
\end{proof}

\subsection{Filtering by Conjugacy Classes of Fixed Sets}\label{filtering}

	Now we complete the proof of the main theorem by showing that
\[ \Omega\mathbb{L}^{-\infty}(\mathcal{B}(CX,CY \cup X,p_X;R))^{h\mathfrak{f}\Gamma} \simeq \Omega\mathbb{L}^{-\infty}(\mathcal{B}(\Sigma X,\Sigma Y,p_X;R))^{h\mathfrak{f}\Gamma}. \]
By Lemma $\ref{lemmaa}$, this is achieved by proving the following.
	
\begin{theo}\label{filtration}
	For every $G \in \mathfrak{f}$,
\[ \mathbb{L}^{-\infty}(\mathcal{B}(CX,CY \cup X,p_X;R))^G \simeq \mathbb{L}^{-\infty}(\mathcal{B}(\Sigma X,\Sigma Y,p_X;R))^G. \]
\end{theo}

	Given $G \in \mathfrak{f}$, consider the subgroup lattice for $G$, and let $H\leq G$. Define the distance from $G$ to $H$, ${\rm dist}(H)$, to be the maximum number of steps needed to reach $H$ from $G$ on the subgroup lattice. Notice that ${\rm dist}(gHg^{-1})={\rm dist}(H)$ for each $g\in G$.

	Let $n={\rm dist}(1)$, and let $l_i$ be the number of conjugacy classes of subgroups with distance $i$ from $G$. For each $i$, $1\leq i\leq n-1$, choose a representative, $H_{i,j}$, $1\leq j\leq l_i$, for each of the conjugacy classes with distance $i$ from $G$. Order the $H_{i,j}$'s using the dictionary order on the indexing set. Now re-index the sequence of subgroups according to the ordering so that
\[H_1=H_{1,1} \ \ , \ \ H_2=H_{1,2} \ \ , \ ... \ , \ \ H_m=H_{n-1,l_{n-1}}. \]
For convenience define $H_0=G$ and $H_{m+1}=1$.

	For each $j$, $0\leq j\leq m+1$, define $C_j=\bigcup_{g \in G} X^{gH_jg^{-1}}$, which will often be referred to as a {\it conjugacy class of fixed sets}.

	For each $k$, $0\leq k\leq m+1$, define $Z_k=\bigcup_{0 \leq j \leq k} C_j$. Also define $Y_k=Z_k \cap Y$. Since it is possible that $Y_k= \emptyset$, we define $C(\emptyset)=\{ 0 \}$ and $\Sigma (\emptyset)=\{ 0,1 \}$ for the convenience of notation.
	
	Notice that $Z_k$ is contractible for every $k$, $0\leq k\leq m+1$. Also notice that since $gX^H=X^{gHg^{-1}}$, each of the conjugacy classes of fixed sets is $G$-invariant. Therefore $Z_k$ is $G$-invariant for every $k$, $0\leq k\leq m+1$.

	Since taking fixed sets commutes with applying $\mathbb{L}^{-\infty}$, we will be working with fixed categories. To simplify the notation set
\begin{align}
	&\mathcal{B}^G_R(C(Z_k))= \mathcal{B}^G(C(Z_k),C(Y_k) \cup Z_k,p_{Z_k};R); \notag\\ 
	&\mathcal{B}^G_R(\Sigma (Z_k)) = \mathcal{B}^G(\Sigma (Z_k),\Sigma (Y_k),p_{Z_k};R); \notag\\
	&\mathcal{B}^G_R(C(Z_k);R)_{Z_{k-1}} = \mathcal{B}^G_R(C(Z_k))_{C(Y_{k-1}) \cup Z_{k-1}}; \notag\\
        &\mathcal{B}^G_R(\Sigma (Z_k))_{Z_{k-1}} = \mathcal{B}^G_R(\Sigma (Z_k))_{\Sigma (Y_{k-1})}; \notag\\
	&\mathcal{B}^G_R(C(Z_k))^{>Z_{k-1}} = \mathcal{B}^G_R(C(Z_k))^{(C(Y_k) \cup Z_k)-(C(Y_{k-1}) \cup Z_{k-1})}; \notag\\
	&\mathcal{B}^G_R(\Sigma (Z_k))^{>Z_{k-1}} = \mathcal{B}^G_R(\Sigma (Z_k))^{\Sigma(Y_k)-\Sigma (Y_{k-1})} . \notag
\end{align}

	Theorem $\ref{filtration}$ is proved by induction on the chain \[X^G=Z_0 \subseteq Z_1 \subseteq \dddot{} \subseteq Z_m \subseteq Z_{m+1}=X. \]
As in~\cite{me2}, this is accomplished with the following three lemmas.

\begin{lem}\label{lemma1}
	Let $1\leq k\leq m+1$. Then there are involution preserving equivalences
	\begin{enumerate}
		\item[i.] $\mathcal{B}^G_R(C(Z_{k}))_{Z_{k-1}} \cong \mathcal{B}^G_R(C(Z_{k-1}))$, and
		\item[ii.] $\mathcal{B}^G_R(\Sigma (Z_{k}))_{Z_{k-1}} \cong \mathcal{B}^G_R(\Sigma (Z_{k-1}))$.
	\end{enumerate}
\end{lem}

\begin{proof}
	This lemma was proved in~\cite{me2} without concern for the involution. The basic idea from the proof still works but needs to be taken a little further.

	Consider part (i). We prove this equivalence by constructing an equivariant function $f$, from $C(Z_k)$ into itself, that is not quite an eventually continuous map but is close enough to induce an involution preserving functor from $\mathcal{B}^G_R(C(Z_{k}))_{Z_{k-1}}$ to $\mathcal{B}^G_R(C(Z_{k-1}))$ that is an inverse to the inclusion functor up to natural equivalence. The only difficulty in defining such a function appears when there are sequences consisting of points outside of $C(Z_{k-1})$ that converge to points in $C(Y_{k-1}) \cup Z_{k-1}$. If this does not happen, then every point not in $C(Z_{k-1})$ is greater than a fixed distance away from $C(Y_{k-1}) \cup Z_{k-1}$. Since objects are zero in a neighborhood of $C(Y_k) \cup Z_k-(C(Y_{k-1}) \cup Z_{k-1})$, there will only be finitely many points in $C(Z_k)-C(Z_{k-1})$ for which an object is non-zero. Thus, we could define $f$ to be the identity on $C(Y_k)\cup Z_k \cup C(Y_{k-1})$ and to send every other point to a chosen point in $\mathcal{E}^G \times (0,1)$.
	
	Assume that there are sequences consisting of points outside of $C(Z_{k-1})$ that converge to points in $C(Y_{k-1}) \cup Z_{k-1}$.. Choose a representative $x$, for each orbit not contained in $C(Z_{k-1})$, such that $G_x=H_k$. Choose $y_x\in (C(Y_{k-1}) \cup Z_{k-1})^{H_k}$ such that $d(x,y_x)=d(x,(C(Y_{k-1}) \cup Z_{k-1})^{H_k})$. This is possible by the assumption that there exist sequences consisting of points outside of $C(Z_{k-1})$ converging to points in $C(Y_{k-1}) \cup Z_{k-1}$. Choose $z_x\in (Z_{k-1}-Y_{k-1})^{H_k}\times (0,1)$ such that $d(z_x,y_x)\leq d(x, C(Y_{k-1}) \cup Z_{k-1})$. This choice is possible by the fourth condition of the main theorem. Define $f:(C(Z_k),C(Y_k\cup Z_k)) \to (C(Z_k),C(Y_k\cup Z_k))$ by
\[ f(a)=\left\{
\begin{array}{cl}
	gz_x & {\rm if} \ a=gx \\
	a & {\rm if} \ a \in C(Y_k)\cup Z_k \cup C(Y_{k-1}).
\end{array}
\right.\]
This is well-defined since $g_1x=g_2x$ if and only if $g_2^{-1}g_1\in H_k$. Thus, $f$ is $G$-equivariant. As mentioned before, $f$ is not eventually continuous. More precisely, the second and third conditions of eventually continuous maps are not satisfied by $f$. The second condition ensures that the induced functor sends objects to objects. To replace this condition, we show that the inverse image under $f$ of a compact set in $C(Z_{k-1})-(C(Y_{k-1}) \cup Z_{k-1})$ contains only finitely many points over which an object can be non-zero. We do this as follows.

	Let $w\in C(Z_{k-1})-(C(Y_{k-1}) \cup Z_{k-1})$ and $d>0$ such that $0<d<d(w,C(Y_{k-1}) \cup Z_{k-1})$ be given. Let
\[ U=\big\{ z\in C(Z_k) \, \big| \, d(z,C(Y_{k-1}) \cup Z_{k-1})<d(w,C(Y_{k-1}) \cup Z_{k-1}-d) \big\}. \]
If $B_d(w)$ is the open ball in $C(Z_{k-1})$ about $w$ of radius $d$, then $f^{-1}(\overline{B_d(w)})\subseteq C(Z_k)-U$. To see why, let $a\in f^{-1}(\overline{B_d(w)})$. If $a\in C(Z_{k-1})$ then $a=f(a)\in \overline{B_d(w)}$. If $a\notin C(Z_{k-1})$ then $a=gx$, where $x$ is the chosen representative, and $d(w,C(Y_{k-1}) \cup Z_{k-1})-d \leq d(gz_x,C(Y_{k-1}) \cup Z_{k-1})=d(z_x, C(Y_{k-1}) \cup Z_{k-1}) \leq d(z_x, C(Y_{k-1}^{H_k}) \cup Z_{k-1}^{H_k}) \leq d(z_x,y_x) \leq d(x, C(Y_{k-1}) \cup Z_{k-1})=d(a, C(Y_{k-1}) \cup Z_{k-1})$. Given an object $H$, there is a neighborhood $V\subseteq C(Z_k)$ of $C(Y_k) \cup Z_k-(C(Y_{k-1}) \cup Z_{k-1})$ on which $H$ is zero. Therefore,
\[ \big\{a\in \overline{f^{-1}(\overline{B_d(w)})} \, \big| \, H_a \neq 0 \big\} \subset (C(Z_k)-U)\cap (C(Z_k)-V). \]
Since $(C(Z_k)-U)\cap (C(Z_k)-V) \subseteq C(Z_k)-(C(Y_k \cup Z_k)$ is compact, this set is finite.

	The third condition on eventually continuous maps is that they are continuous on the boundary. This condition guarantees that the image under $f_{\sharp}$ of a morphism is continuously controlled.  Since objects are zero on a neighborhood of $C(Y_k) \cup Z_k-(C(Y_{k-1}) \cup Z_{k-1})$, it suffices to show that $f$ is continuous on $C(Y_{k-1}) \cup Z_{k-1}$. Let $\{ x_n \}$ be a sequence in $C(Z_k)-C(Z_{k-1})$ converging to $y\in C(Y_{k-1}) \cup Z_{k-1}$. Since $G$ is a finite group, all but finitely many terms in the sequence must be contained in a subsequence of the form $\{ gx_m \}$ for a fixed $g\in G$, where each $x_k$ is the chosen representative of its orbit. Since $\{ x_m \}$ converges to $g^{-1}y$ and $G_{x_m}=H_k$, $g^{-1}y$ is fixed by $H_k$. Thus, $\{ y_{x_m} \}$ converges to $g^{-1}y$. Therefore, $\{ z_{x_m} \}$ converges to $g^{-1}y$ and $\{ f(gx_m)\}=\{ gz_{x_m} \}$ converges to $y=f(y)$. Since $f|_{C(Y_k) \cup Z_k}=1_{C(Y_k) \cup Z_k}$, $f$ induces a functor that is an inverse to the inclusion functor up to natural equivalence, as with eventually continuous maps.
	
	This completes the proof of part (i). Part (ii) is proven the same way replacing cones with suspensions.
\end{proof}

\begin{lem}\label{lemma2}
	Let $1\leq k\leq m+1$. Then there are involution preserving equivalences
	\begin{enumerate}
		\item[i.] $\mathcal{B}^G_{R}(C(Z_{k}))^{>Z_{k-1}} \cong \mathcal{B}_{RH_k}(C(Z_{k} / G))^{>Z_{k-1} / G}$, and
		\item[ii.] $\mathcal{B}^G_{R}(\Sigma (Z_{k}))^{>Z_{k-1}} \cong \mathcal{B}_{RH_k}(\Sigma (Z_{k} / G))^{>Z_{k-1} / G}$.
	\end{enumerate}
\end{lem}

\begin{proof}
	These equivalences are proven in~\cite{me2}. However, we need to check that the involution is preserved. We begin by recalling the argument from the proof of~\cite[Lemma 7.4]{me2}. Since we are taking germs away from $C(Y_{k-1}) \cup Z_{k-1}$ and $C(Y_{k-1} / G) \cup Z_{k-1} / G$, every morphism has a representative that is zero on $C(Z_{k-1})$ and on $C(Z_{k-1} / G)$, respectively. It is therefore irrelevant what the objects over $C(Z_{k-1})$ and $C(Z_{k-1} / G)$ are. Also, the part of an object over a point not in $C(Y_{k-1}) \cup Z_{k-1}$ must be an $R[gH_kg^{-1}]$-module, for some $g\in G$. This explains how the objects in $\mathcal{B}^G_{R}(C(Z_{k}))^{>Z_{k-1}}$ and $\mathcal{B}_{RH_k}(C(Z_{k} / G))^{>Z_{k-1} / G}$ are identified. Since we are taking germs, the components of a morphism need to become small. Therefore non-zero components of a morphism have the same isotropy, namely a conjugate of $H_k$. Furthermore, the equivariance of morphisms in $\mathcal{B}^G_{R}(C(Z_{k}))^{>Z_{k-1}}$ implies that there is only one choice when lifting a morphism from $\mathcal{B}_{RH_k}(C(Z_{k} / G))^{>Z_{k-1} / G}$. In other words, up to equivariance, the components of a morphism in $\mathcal{B}^G_{R}(C(Z_{k}))^{>Z_{k-1}}$ are the same as the morphism that it is identified with in $\mathcal{B}_{RH_k}(C(Z_{k} / G))^{>Z_{k-1} / G}$. But this implies that conjugate transposes are sent to conjugate transposes. Therefore the involution is preserved. This proves the first part of this lemma. The same argument proves the second part, replacing cones with suspensions.
\end{proof}

\begin{lem}\label{lemma3}
	For each $k$, $1\leq k\leq m+1$,
\[ \mathbb{L}^{-\infty}\big(\mathcal{B}_{RH_k}(C(Z_{k} / G))^{>Z_{k-1} / G}\big) \simeq \mathbb{L}^{-\infty}\big(\mathcal{B}_{RH_k}(\Sigma (Z_{k} / G))^{>Z_{k-1} / G}\big). \]
\end{lem}

\begin{proof}
	The proof proceeds just as in~\cite[Lemma 7.5]{me2}. Consider the following commutative diagram:
\[ \xymatrix{
	\mathcal{B}_{RH_k}(C({Z_k} / G))_{Z_{k-1}} \ar[d]_a \ar[r] & \mathcal{B}_{RH_k}(C({Z_k} / G)) \ar[d]_b \ar[r] & \mathcal{B}_{RH_k}(C({Z_k} / G))^{>Z_{k-1}} \ar[d]_c \\
	\mathcal{B}_{RH_k}(\Sigma({Z_k} / {G}))_{Z_{k-1}} \ar[r] & \mathcal{B}_{RH_k}(\Sigma({Z_k} / {G})) \ar[r] & \mathcal{B}_{RH_k}(\Sigma({Z_k} / {G}))^{>Z_{k-1}}. } \]
Since $Z_k$ is contractible, $b$ induces a weak homotopy equivalence by Corollary~\ref{zero} and Proposition~\ref{key}. By Lemma~\ref{lemma1}, we can apply Corollary~\ref{zero} and Proposition~\ref{key} again to see that $a$ induces a weak homotopy equivalence. Since each row in the diagram is a Karoubi filtration, the Five Lemma shows that $c$ also induces a weak homotopy equivalence.
\end{proof}

We are now ready to prove Theorem~\ref{filtration}.

\begin{proof}[Proof of Theorem~\ref{filtration}]
	Since $G$ acts trivially on $X^G$, $\mathcal{B}^G_R(C(X^G)) \cong \mathcal{B}_{RG}(C(X^G))$ and $\mathcal{B}^G_R(\Sigma(X^G)) \cong \mathcal{B}_{RG}(\Sigma(X^G))$. Since $X^G$ is contractible, Proposition~\ref{key} implies ${\mathbb{L}}^{-\infty}(\mathcal{B}_{RG}(C(X^G))) \simeq {\mathbb{L}}^{-\infty}(\mathcal{B}_{RG}(\Sigma(X^G)))$. Therefore
\[ {\mathbb{L}}^{-\infty}\big(\mathcal{B}^G_R(C(X^G))\big) \simeq \mathbb{L}^{-\infty}\big(\mathcal{B}^G_R(\Sigma(X^G))\big). \]
This completes the base case of the induction. 
	
	Assume now that $\mathbb{L}^{-\infty}(\mathcal{B}^G_R(C(Z_{k-1}))) \simeq \mathbb{L}^{-\infty}(\mathcal{B}^G_R(\Sigma(Z_{k-1})))$. We want to show that $\mathbb{L}^{-\infty}(\mathcal{B}^G_R(C(Z_k))) \simeq \mathbb{L}^{-\infty}(\mathcal{B}^G_R(\Sigma(Z_k)))$.
Consider the following commutative diagram:
\[ \xymatrix{
	\mathcal{B}^G_R(C(Z_k))_{Z_{k-1}} \ar[d]_a \ar[r] & \mathcal{B}^G_R(C(Z_k)) \ar[d]_b \ar[r] & \mathcal{B}^G_R(C(Z_k))^{>Z_{k-1}} \ar[d]_c \\
	\mathcal{B}^G_R(\Sigma(Z_k))_{Z_{k-1}} \ar[r] & \mathcal{B}^G_R(\Sigma(Z_k)) \ar[r] & \mathcal{B}^G_R(\Sigma(Z_k))^{>Z_{k-1}}. } \]

	By Lemma $\ref{lemma1}$, Proposition $\ref{key}$, and the induction hypothesis, $a$ induces a weak homotopy equivalence. By Lemmas $\ref{lemma2}$ and $\ref{lemma3}$, $c$ induces a weak homotopy equivalence. Since each row in the diagram is a Karoubi filtration, the Five Lemma shows that $b$ also induces a weak homotopy equivalence.
 \end{proof}

	This completes the proof of the main theorem.

\end{document}